\documentclass[10pt]{article}
\usepackage{graphicx}
\usepackage{amsmath,amssymb,amsthm,amsfonts}
\usepackage{amssymb}
\usepackage{mathtools}
\newtheorem{thm}{Theorem}[section]

\newtheorem{lemma}[thm]{Lemma}

\theoremstyle{definition}
\newtheorem{defn}{Definition}[section]
\theoremstyle{remark}
\newtheorem{rem}{Remark}[section]

\numberwithin{equation}{section}

\DeclareMathSymbol{\C}{\mathalpha}{AMSb}{"43}

\newcommand{\bsub}{\begin{subequations}}
\newcommand{\esub}{\end{subequations}$\!$}

\begin{document}

\title{A variational model with fractional-order regularization term arising in registration of diffusion tensor image \thanks{Email: hanhuan11@mails.ucas.ac.cn (H. Han)
}\thanks{This work was supported by NSFC under grant No.11471331 and partially supported by National Center for Mathematics and Interdisciplinary Sciences.}}

\author{Huan Han\\
\small Wuhan Insitute of Physics and Mathematics, Chinese Academy of Sciences,\\ \small  P.O. Box 71010, Wuhan, 430071, China
}

\date{}

\smallbreak \maketitle

\begin{abstract} In this paper, a new variational model with fractional-order regularization term arising in registration of diffusion tensor image(DTI) is presented. Moreover, the existence of its solution is proved to ensure that there is a regular solution for this model.
\end{abstract}

\vskip 0.2truein



\vskip 0.2truein
\section{Introduction}
Let $\Omega\subset \mathbb{R}^3$ be an open bounded domain, i.e., $\Omega=(a_1,b_1)\times(a_2,b_2)\times(a_3,b_3)$. Suppose $T$ and $D$ are two functions defined from $\Omega$ to the set of $3\times 3$ real Symmetric Positive Definite matrixes($SPD(3)$ in short). That is,
\begin{align}\label{eq1-1}
T: \Omega\rightarrow SPD(3),\ \ D: \Omega\rightarrow SPD(3).
\end{align}

In DTI registration, $T$ and $D$ are viewed as two images defined on $\Omega$, where $T$ is called floating image and $D$ is called target image. The goal of registration is to find a $1$-to-$1$ spatial transformation $h:\Omega\rightarrow\Omega$ such that $T\circ h(\cdot)$ is close to $D(\cdot)$ in some sense. On the other hand, in order to keep $T\circ h(\cdot)$ align with spatial transformation, reorientation of $T\circ h(\cdot)$ must be additionally considered. For this purpose, Alexander\cite{APBG} put forward two reorientation strategies: finite strain(FS) strategy and preservation principle direction(PPD) strategy. Based on FS strategy, Li\cite{LSTDWT} introduced a new operator ``$\diamond$'' defined by
\begin{align}\label{eq1-2real}
T\diamond h(x)=R[T\circ h(x)]R^T\ \mathrm{with}\ R=J^T(JJ^T)^{-\frac{1}{2}}\ \mathrm{and}\ J=\nabla_x h^{-1}(x).
\end{align}

With the help of this operator, the DTI registration model(cf. \cite{HH}) can be formulated as
\begin{align}\label{eq1-2}
\bar{{v}}=\arg\min_{v} \bar{H}(v),
\end{align}
where $\bar{H}(v)=\int_0^\tau\|Lv(\cdot,t)\|_{L^2(\Omega)}^2 dt+\|T\diamond h(\cdot)-D(\cdot)\|_{L^2(\Omega)}^2$, $L:[H^3_0(\Omega)]^3\rightarrow [L^2(\Omega)]^3$ is a linear differential operator satisfing
\begin{align}\label{eqnew1}
&\|Lv(\cdot,t)\|^2_{L^2(\Omega)}\triangleq\sum_{i=1}^3\int_{\Omega}|(Lv)_i(x,t)|^2dx\geq c \|v(\cdot,t)\|^2_{[H^3_0(\Omega)]^3},&
\end{align}
for some constant $c>0$, and
\begin{align} \label{eq1-3}
&\frac{d{\eta}(s;t,x)}{ds}=v(\eta(s;t,x),s),\ \eta(t;t,x)=x\ \ \mathrm{and } \ \ h(x)=\eta(0;\tau,x).&
\end{align}
\begin{rem}
Note that in (\ref{eq1-2real}), $R$, $J$ and $T\circ h(x)$ are all $3\times 3$ matrixes with $h^{-1}(x)=\eta(\tau;0,x)$ and $J=\nabla_x\eta(\tau;0,x)$ ({cf. \cite{HH}}). Here the existence of function $h^{-1}$ is given by {\bf(ii)} in Lemma \ref{lemm3-5} and the definition of $(JJ^T)^{-\frac{1}{2}}$ can refer to Appendix in \cite{HH}.\qed
\end{rem}

In \cite{HH}, authors prove that there exists a solution to variational problem (\ref{eq1-2})-(\ref{eq1-3}) on some suitable space.
Note that almost all the DTI registration model\cite{LSTDWT,YVFPPGAC} have employed integer-order derivatives in linear  differential operator $L$. In fact, during the last decades, it has been showed that many problems involving science and engineering can be modeled more accurately by employing fractional-order derivatives\cite{PI,WD,ZC} than integer-order derivatives. Motivated by this fact, the aim of this paper is to employ fractional-order derivatives in DTI registration model.

Before giving our results, we introduce some notations and definitions.

Throughout this paper, we define $\Omega\triangleq(a_1,b_1)\times(a_2,b_2)\times(a_3,b_3)\subset \mathbb{R}^3$. Moreover, for $x\in\Omega$, the inner product and modulus of matrix $A(x)=\left(a_{ij}(x)\right)_{n\times m}$, $B=\left(b_{ij}(x)\right)_{n\times m}$ are defined as
\begin{align}
A(x)\cdot B(x)=\sum_{i=1}^n\sum_{j=1}^ma_{ij}(x)b_{ij}(x),\ \ \|A(x)\|=\sqrt{\sum_{i=1}^n\sum_{j=1}^m a_{ij}^2(x)},\nonumber
\end{align}
respectively.

Furthermore, we say $A(x)$ is continuous on $\Omega$ if $a_{ij}(x)(i=1,2,\cdots,n;j=1,2,\cdots,m)$ are continuous on $\Omega$.

Moreover, for matrix sequence $A_{k}(x)=\left(a_{ij}^k(x)\right)_{n\times m}$, we say $A_{k}(x)\xrightarrow{k} A(x)$ if $a_{ij}^k(x)\xrightarrow{k} a_{ij}(x)(i=1,2,\cdots,n;j=1,2,\cdots,m)$.

Based on definition of Riemann-Liouville derivative in \cite{ER}, for $x=(x_1,x_2,x_3)\in\Omega$ and function $f:\Omega\rightarrow \mathbb{R}$, define

\begin{align}
\frac{\partial^\alpha f(x)}{\partial x_i^\alpha}\triangleq D^\alpha_{[a_i,x_i]}f(x)=\frac{1}{\Gamma([\alpha]+1-\alpha)}\left(\frac{d}{dx_i}\right)^{[\alpha]+1}\int_{a_i}^{x_i}\frac{f^{(i)}(x,t)}{(x_i-t)^{\alpha-[\alpha]}}dt,\nonumber
\end{align}
\begin{align}
\frac{\partial^{\alpha*} f(x)}{\partial x_i^{\alpha*}}\triangleq D^{\alpha}_{[x_i,b_i]}f(x)=\frac{1}{\Gamma([\alpha]+1-\alpha)}\left(-\frac{d}{dx_i}\right)^{[\alpha]+1}\int_{x_i}^{b_i}\frac{f^{(i)}(x,t)}{(t-x_i)^{\alpha-[\alpha]}}dt,\nonumber
\end{align}
where $\Gamma(s)=\int_0^{+\infty}x^{s-1}e^{-x}dx$, $[\cdot]$ is round down function, here and in what follows, $f^{(1)}(x,t)=f(t,x_2,x_3)$, $f^{(2)}(x,t)=f(x_1,t,x_3)$, $f^{(3)}(x,t)=f(x_1,x_2,t)$  and $i=1,2,3$.

\begin{defn}\label{def1-2}
For $\alpha>0$ and function $g:\Omega\rightarrow \mathbb{R}$, define semi-norms
\begin{align}
|g|_{F^\alpha_L(\Omega)}=\left(\int_{\Omega}\|\nabla^\alpha g(x)\|^2 dx\right)^{\frac{1}{2}}, |g|_{F^\alpha_R(\Omega)}=\left(\int_{\Omega}\|\nabla^{\alpha *} g(x)\|^2 dx\right)^{\frac{1}{2}},\nonumber
\end{align}
and norms
\begin{align}
\|g\|_{F^\alpha_L(\Omega)}=\left(\|g\|_{L^2(\Omega)}^2+|g|^2_{F^\alpha_L(\Omega)}\right)^{\frac{1}{2}},\|g\|_{F^\alpha_R(\Omega)}=\left(\|g\|_{L^2(\Omega)}^2+|g|^2_{F^\alpha_R(\Omega)}\right)^{\frac{1}{2}},\nonumber
\end{align}
where $\nabla^\alpha g(x)=\left(\frac{\partial^\alpha g(x)}{\partial x_i^\alpha}\right)_{1\times 3}$ and  $\nabla^{\alpha *} g(x)=\left(\frac{\partial^{\alpha*} g(x)}{\partial x_i^{\alpha*}}\right)_{1\times 3}$.
\end{defn}

Based on Definition \ref{def1-2}, define space $F^\alpha_{L,0}(\Omega)$ and $F^\alpha_{R,0}(\Omega)$ as the closure of $C_0^{\infty}(\Omega)$ under the norm $\|\cdot\|_{F^\alpha_{L,0}(\Omega)}$ and $\|\cdot\|_{F^\alpha_{R,0}(\Omega)}$, respectively.

\begin{defn}\label{1def3}
 For $\alpha>0$ and $u\in L^1(\mathbb{R}^3)$, define the semi-norm and norm
\begin{align}
|u|_{H^\alpha(\mathbb{R}^3)}=\left(\int_{\mathbb{R}^3}\|\xi\|^{2\alpha}|\widehat{{u}}(\xi)|^2d\xi\right)^{\frac{1}{2}},\ \ \ \ \|u\|_{H^\alpha(\mathbb{R}^3)}=\left(\|u\|^2_{L^2(\mathbb{R}^3)}+|u|^2_{H^\alpha(\mathbb{R}^3)}\right)^{\frac{1}{2}},\nonumber
\end{align}
where here and in what follows, $\widehat{{u}}(\xi)=\frac{1}{(2\pi)^{\frac{3}{2}}}\int_{\mathbb{R}^3}u(x)e^{-x\cdot\xi}dx$.

Define Sobolev space $H^\alpha(\mathbb{R}^3)$  as the closure of $C^{\infty}_0(\mathbb{R}^3)$ under the norm $\|\cdot\|_{H^\alpha(\mathbb{R}^3)}.$
\end{defn}

In Definition \ref{1def3}, if we restrict $\mathbb{R}^3$ to $\Omega$, then $H^\alpha_0(\Omega)$ is the closure of $C_0^{\infty}(\Omega)$ under the norm $\|\cdot\|_{H^\alpha_0(\Omega)}$.

For $\tau>0$, $\alpha>2.5$, $\alpha\neq m+0.5$, $m\in \mathbb{N}$ and $u:\Omega\times[0,\tau]\rightarrow \mathbb{R}^3$, define a separable Hilbert space
\begin{align}
\mathcal{F}\triangleq\{u(x,t)=(u_i(x,t))_{1\times 3}:u_i(x,t)\in F^\alpha_{L,0}(\Omega)\  \mathrm{for\ any} \ t\in[0,\tau] \ \mathrm{and}\ i=1,2,3\},\nonumber
\end{align}
endowing with the following inner product and norm
\begin{align}
&(u,v)_{{\mathcal{F}}}=\int_0^\tau \int_\Omega \nabla^\alpha u(x,t)\cdot \nabla^\alpha v(x,t)dxdt, \ \ \ \|u\|_{\mathcal{F}}^2=\int_0^\tau\|\nabla^\alpha u(\cdot,t)\|^2_{L^2(\Omega)}dt,\nonumber&
\end{align}
where $\nabla^\alpha u(x,t)=\left(\frac{\partial^\alpha u_i(x,t)}{\partial x^\alpha_j}\right)_{3\times 3}$.

Based on the above notations and definitions, the variational model with fractional-order regularization term arising in registration of diffusion tensor image(DTI) can be formulated as
\begin{align}\label{eq1-4}
\bar{{v}}=\arg\min_{v\in\mathcal{F}} {H}(v),
\end{align}
where ${H}(v)=\int_0^\tau\|\nabla^\alpha v(\cdot,t)\|_{L^2(\Omega)}^2 dt+\|T\diamond h(\cdot)-D(\cdot)\|_{L^2(\Omega)}^2$ and $h(x)$ is defined by (\ref{eq1-3}).

Another purpose of this paper is to give a rigid proof on the existence of solution to (\ref{eq1-4}). As to this problem, we have the following result:
\begin{thm}\label{th-A10}
Let $T$ and $D$ be two functions defined by (\ref{eq1-1}), and let the set $\triangle_T\triangleq\{x:T(\cdot) \ \ \mathrm{ is \ discontinuous\ at} \ x\}$ be a set of measure zero. If $\mathop {\max }\limits_{x\in \Omega} \|T(x)\|<+\infty$, $G\triangleq\mathop {\max }\limits_{x\in \Omega} \|T(x)-D(x)\|^2<+\infty$, then the variational problem (\ref{eq1-4}) admits a solution $\bar{{v}}(x,s)\in {\mathcal{F}}$ with $\|\bar{{v}}\|_{\mathcal{F}}^2\leq G|\Omega|$. Furthermore,  by (\ref{eq1-3}), $\bar{{v}}(x,s)$ induces a $1$-to-$1$ and onto mapping $\bar{{h}}(x)\in [C^{[\alpha-1.5],\lambda}(\Omega)]^3$ defined from $\Omega$ to $\Omega$, where $0<\lambda\leq\alpha-1.5-[\alpha-1.5]$. Moreover, $\nabla_x\bar{{h}}(x)$ is given by Lemma \ref{lemm3-5}.
\end{thm}

\begin{rem}
In fact, if $H(v)$ in (\ref{eq1-4}) is formulated as
\begin{align}
{H}(v)=\int_0^\tau\|\nabla^{\alpha*} v(\cdot,t)\|_{L^2(\Omega)}^2 dt+\|T\diamond h(\cdot)-D(\cdot)\|_{L^2(\Omega)}^2,\nonumber
\end{align}
then there also exists a global minimizer to $H(v)$ on space
\begin{align}
\mathcal{F}_1\triangleq\{u(x,t)=(u_i(x,t))_{1\times 3}:u_i(x,t)\in F^\alpha_{R,0}(\Omega)\  \mathrm{for\ any} \ t\in[0,\tau]\  \mathrm{and}\ i=1,2,3\},\nonumber
\end{align}
endowing with the following inner product and norm
\begin{align}
&(u,v)_{{\mathcal{F}}_1}=\int_0^\tau \int_\Omega \nabla^{\alpha*} u(x,t)\cdot \nabla^{\alpha*} v(x,t)dxdt, \ \ \  \|u\|_{\mathcal{F}_1}^2=\int_0^\tau\|\nabla^{\alpha*} u(\cdot,t)\|^2_{L^2(\Omega)}dt,\nonumber&
\end{align}
where $\alpha>2.5$, $\alpha\neq m+0.5$, $m\in \mathbb{N}$ and $\nabla^{\alpha*} u(x,t)=\left(\frac{\partial^{\alpha*} u_i(x,t)}{\partial x^{\alpha*}_j}\right)_{3\times 3}$.
\end{rem}

\section{Equivalence of $F^\alpha_{L,0}(\Omega)$, $F^\alpha_{R,0}(\Omega)$ and $H^\alpha_0(\Omega)$}

In \cite{HH}, authors impose the condition (\ref{eqnew1}) on $L$ such that $v(\cdot,t)\in [H^3_0(\Omega)]^3\hookrightarrow [C^1(\Omega)]^3$ which ensures the existence and uniqueness of solution to (\ref{eq1-3}). As the basic space of this paper, $F^\alpha_{L,0}(\Omega)$ and $F^\alpha_{R,0}(\Omega)$ are also needed to embedded into $C^1(\Omega)$. Otherwise, the uniqueness of solution to (\ref{eq1-3}) can not be guaranteed\cite{T}.

For this purpose, we will prove the equivalence of $F^\alpha_{L,0}(\Omega)$, $F^\alpha_{R,0}(\Omega)$ and $H^\alpha_0(\Omega)$, since $H^\alpha(\Omega)\hookrightarrow C^1(\Omega)(\alpha>2.5)$(cf. [2, Theorem 4.57]).

First, we introduce some definitions.

\begin{defn}
For $\alpha>0$ and function $g:\mathbb{R}^3\rightarrow \mathbb{R}$, define the semi-norms
\begin{align}
|g|_{F^\alpha_L(\mathbb{R}^3)}=\left(\int_{\mathbb{R}^3}\|D^\alpha g(x)\|^2dx\right)^{\frac{1}{2}},\ \ |g|_{F^\alpha_R(\mathbb{R}^3)}=\left(\int_{\mathbb{R}^3}\|D^{\alpha *}g(x)\|^2dx\right)^{\frac{1}{2}},\nonumber
\end{align}
and norms
\begin{align}
\|g\|_{F^\alpha_L(\mathbb{R}^3)}=\left(\|g\|_{L^2(\mathbb{R}^3)}^2+|g|^2_{F^\alpha_L(\mathbb{R}^3)}\right)^{\frac{1}{2}},\|g\|_{F^\alpha_R(\mathbb{R}^3)}=\left(\|g\|_{L^2(\mathbb{R}^3)}^2+|g|^2_{F^\alpha_R(\mathbb{R}^3)}\right)^{\frac{1}{2}},\nonumber
\end{align}
where $D^\alpha g(x)=(D^\alpha_j g(x))_{1\times 3}$, $D^{\alpha*} g(x)=(D^{\alpha*}_j g(x))_{1\times 3}$ and
\begin{align}
D^\alpha_{j}g(x)=\frac{1}{\Gamma([\alpha]+1-\alpha)}\left(\frac{d}{dx_j}\right)^{[\alpha]+1}\int_{-\infty}^{x_j}\frac{g^{(j)}(x,t)}{(x_j-t)^{\alpha-[\alpha]}}dt,\\
D^{\alpha*}_{j}g(x)=\frac{1}{\Gamma([\alpha]+1-\alpha)}\left(-\frac{d}{dx_j}\right)^{[\alpha]+1}\int_{x_j}^{+\infty}\frac{g^{(j)}(x,t)}{(t-x_j)^{\alpha-[\alpha]}}dt.
\end{align}

Define $F^\alpha_L(\mathbb{R}^3)$, $F^\alpha_R(\mathbb{R}^3)$ as the closure of $C^\infty_0(\mathbb{R}^3)$ under the norm $\|\cdot\|_{F^\alpha_L(\mathbb{R}^3)}$ and $\|\cdot\|_{F^\alpha_R(\mathbb{R}^3)}$, respectively.
\end{defn}

\begin{defn}
For $\alpha>0$ and function $g:\mathbb{R}^3\rightarrow \mathbb{R}$, define semi-norm
\begin{align}\label{eq2-1}
|g|_{F^\alpha_S(\mathbb{R}^3)}=\left|\int_{\mathbb{R}^3}D^\alpha g(x)\cdot D^{\alpha*} g(x)dx\right|^{\frac{1}{2}},
\end{align}
and norm
\begin{align}\label{eq2-2}
\|g\|_{F^\alpha_S(\mathbb{R}^3)}=\left(\|g\|_{L^2(\mathbb{R}^3)}^2+|g|^2_{F^\alpha_S(\mathbb{R}^3)}\right)^{\frac{1}{2}}.
\end{align}

Define $F^\alpha_S(\mathbb{R}^3)$ as the closure of $C^\infty_0(\mathbb{R}^3)$ under the norm $\|\cdot\|_{F^\alpha_S(\mathbb{R}^3)}$.

If we restrict $\mathbb{R}^3$ to $\Omega$ and replace $D^\alpha g(x)$, $D^{\alpha *}g(x)$ with $\nabla^\alpha g(x)$, $\nabla^{\alpha *}g(x)$ in (\ref{eq2-1}) respectively, then $F^\alpha_{S,0}(\Omega)$ is defined similarly.

\end{defn}

\begin{lemma}\label{lem2-1}
For $\alpha>0$ and $\xi=(\xi_1,\xi_2,\xi_3)\in \mathbb{R}^3$, there holds
\begin{align}\label{eq2-new}
\frac{1}{3}\|\xi^\alpha\|^2\leq\|\xi\|^{2\alpha}\leq3\|\xi^\alpha\|^2,
\end{align}
where $\xi^\alpha=(\xi_1^\alpha,\xi_2^\alpha,\xi_3^\alpha)$.
\end{lemma}
{\bf Proof.} Since $\|\xi\|^{2\alpha}=\left(\xi_1^2+\xi_2^2+\xi_3^2\right)^\alpha\geq \xi_j^{2\alpha}(j=1,2,3)$, then there holds
\begin{align}
\frac{1}{3}\|\xi^\alpha\|^2=\frac{1}{3}(\xi_1^{2\alpha}+\xi_2^{2\alpha}+\xi_3^{2\alpha})\leq\|\xi\|^{2\alpha}.
\end{align}

On the other hand, $\|\xi\|^{2\alpha}\leq 3\max\limits_{j=1,2,3}\xi_j^{2\alpha}$ implies that
\begin{align}
\|\xi\|^{2\alpha}\leq 3(\xi_1^{2\alpha}+\xi_2^{2\alpha}+\xi_3^{2\alpha})=3\|\xi^\alpha\|^2.
\end{align}

\begin{lemma}\label{lem2-2}
Assume $\alpha>0$, $\xi=(\xi_1,\xi_2,\xi_3)\in \mathbb{R}^3$, then
\begin{align}
\widehat{(D^\alpha_j f)}(\xi)=(i\xi_j)^\alpha \widehat{{f}}(\xi),\ \ \widehat{(D^{\alpha*}_j f)}(\xi)=(-i\xi_j)^\alpha \widehat{{f}}(\xi)\ (j=1,2,3).
\end{align}
\end{lemma}
{\bf Proof.} For function $g:\mathbb{R}\rightarrow \mathbb{R}$ and $\chi\in\mathbb{R}$, by Appendix in \cite{ER}, we know
\begin{align}\label{eq2-9}
\widehat{(D^\alpha g)}(\chi)=(i\chi)^\alpha\widehat{{g}}(\chi),\ \ \widehat{(D^{\alpha *} g)}(\chi)=(-i\chi)^\alpha\widehat{{g}}(\chi),
\end{align}
where $\widehat{{g}}(\chi)=\frac{1}{(2\pi)^{\frac{1}{2}}}\int_{-\infty}^{+\infty}g(y)e^{-iy\chi}dy$.

Based on this conclusion, we have
\begin{align}
\widehat{(D_1^\alpha f)}(\xi)=&\frac{1}{(2\pi)^{\frac{3}{2}}}\int_{-\infty}^{+\infty}\int_{-\infty}^{+\infty}\int_{-\infty}^{+\infty}D_1^\alpha f(x_1,x_2,x_3) e^{-i\sum\limits_{j=1}^3 x_j\xi_j}dx_1dx_2dx_3\nonumber&
\end{align}
\begin{align}
=&\frac{1}{2\pi}\int_{-\infty}^{+\infty}\int_{-\infty}^{+\infty}\left[\frac{1}{(2\pi)^{\frac{1}{2}}}\int_{-\infty}^{+\infty}D_1^\alpha f(x_1,x_2,x_3) e^{-i x_1\xi_1}dx_1\right]e^{-i\sum\limits_{j=2}^3 x_j\xi_j}dx_2dx_3\nonumber&\\
=&\frac{1}{(2\pi)^{\frac{3}{2}}}\int_{-\infty}^{+\infty}\int_{-\infty}^{+\infty}\left[(i\xi_1)^\alpha\int_{-\infty}^{+\infty}f(x_1,x_2,x_3) e^{-i x_1\xi_1}dx_1\right]e^{-i\sum\limits_{j=2}^3 x_j\xi_j}dx_2dx_3\nonumber&\\
=&(i\xi_1)^\alpha\frac{1}{(2\pi)^{\frac{3}{2}}}\int_{-\infty}^{+\infty}\int_{-\infty}^{+\infty}\int_{-\infty}^{+\infty}f(x_1,x_2,x_3) e^{-i\sum\limits_{j=1}^3 x_j\xi_j}dx_1dx_2dx_3\nonumber&\\
=&(i\xi_1)^\alpha \widehat{{f}}(\xi).&
\end{align}

Similarly, we can prove that $\widehat{(D^\alpha_j f)}(\xi)=(i\xi_j)^\alpha \widehat{{f}}(\xi)(j=2,3)$.

With the help of  the second equality of (\ref{eq2-9}), we can prove that
\begin{align}
\widehat{(D^{\alpha*}_j f)}(\xi)=(-i\xi_j)^\alpha \widehat{{f}}(\xi)(j=1,2,3).
\end{align}

\begin{lemma}\label{lem2-3}
Assume $\alpha>0$, then $F^\alpha_L(\mathbb{R}^3)$, $F^\alpha_R(\mathbb{R}^3)$ and $H^\alpha(\mathbb{R}^3)$ are equivalent.
\end{lemma}
{\bf Proof.}  By Lemma \ref{lem2-2} and Plancherel Theorem[5, Theorem 1 in Section 4.3],
\begin{align}
|f|^2_{F^\alpha_L(\mathbb{R}^3)}=&\sum_{j=1}^3\|D_j^\alpha f\|^2_{L^2(\mathbb{R}^3)}
=\sum_{j=1}^3\|\widehat{(D_j^\alpha f)}\|^2_{L^2(\mathbb{R}^3)}
=\sum_{j=1}^3\||\xi_j|^\alpha \widehat{{f}}(\xi)\|^2_{L^2(\mathbb{R}^3)}\nonumber&\\
=&\int_{\mathbb{R}^3}\left|\widehat{{f}}(\xi)\right|^2\sum_{j=1}^3|\xi_j|^{2\alpha}d\xi=\int_{\mathbb{R}^3}\left|\widehat{{f}}(\xi)\right|^2\|\xi^\alpha\|^{2}d\xi.&
\end{align}

By Lemma \ref{lem2-1}, we know that $F^\alpha_L(\mathbb{R}^3)$ and $H^\alpha(\mathbb{R}^3)$ are equivalent.

In a similar way, we can prove that $F^\alpha_R(\mathbb{R}^3)$ and $H^\alpha(\mathbb{R}^3)$ are equivalent.

\begin{lemma}\label{newlem}
Assume $\alpha>0$ and let $f(x)$ be a function defined from $\mathbb{R}^3$ to $\mathbb{R}$, then
\begin{align}\label{eq-newlem}
\int_{\mathbb{R}^3} D^\alpha f(x)\cdot D^{\alpha *}f(x)dx=\cos(\pi\alpha)\|D^\alpha f\|^2_{L^2(\mathbb{R}^3)}.
\end{align}
\end{lemma}
{\bf Proof.} By Parseval equality[5, Theorem 2 in Section 4.3]
\begin{align}
\int_\mathbb{R} u\bar{v}dx=\int_\mathbb{R}\widehat{u}(\xi)\overline{\widehat{v}(\xi)}d\xi,
\end{align}
we have,
\begin{align}\label{hjkljjnn}
\int_{\mathbb{R}^3} D^\alpha_1 f(x)\cdot D^{\alpha *}_1f(x)dx=&\int_{-\infty}^{+\infty}\int_{-\infty}^{+\infty}\int_{-\infty}^{+\infty}D^\alpha_1 f(x) \overline{D^{\alpha *}_1f(x)}dx_1dx_2dx_3\nonumber&\\
=&\int_{-\infty}^{+\infty}\int_{-\infty}^{+\infty}\left[\int_{-\infty}^{+\infty}D^\alpha_1 f(x) \overline{D^{\alpha *}_1f(x)}dx_1\right]dx_2dx_3\nonumber&\\
=&\int_{-\infty}^{+\infty}\int_{-\infty}^{+\infty}\left[\int_{-\infty}^{+\infty} \widehat{(D^\alpha_1 f)}(\xi_1) \overline{\widehat{(D^{\alpha *}_1f)}(\xi_1)}d\xi_1\right]dx_2dx_3,&
\end{align}
where $\widehat{(D^\alpha_1 f)}(\xi_1)=\frac{1}{(2\pi)^{\frac{1}{2}}}\int_{-\infty}^{+\infty}D^\alpha_1 f(x)e^{-ix_1\xi_1}dx_1$, $\widehat{(D^{\alpha*}_1 f)}(\xi_1)=\frac{1}{(2\pi)^{\frac{1}{2}}}\int_{-\infty}^{+\infty}D^{\alpha*}_1 f(x)e^{-ix_1\xi_1}dx_1$.

On the other hand, we know
\begin{equation}\label{eq-216}\overline{(iw)^\alpha}=\left\{\begin{aligned}
  &e^{-i\pi\alpha}\overline{(-iw)^\alpha} \ \ if \ w\geq 0,&\\
  &e^{i\pi\alpha}\overline{(-iw)^\alpha} \ \ \ if \ w<0. &\\
  \end{aligned}\right.\end{equation}

It follows from (\ref{eq-216}) that
\begin{align}
\int_{-\infty}^{+\infty}\widehat{(D^\alpha_1 f)}(\xi_1) \overline{\widehat{(D^{\alpha *}_1f)}(\xi_1)}d\xi_1=&\int_{-\infty}^0 (i\xi_1)^\alpha\widehat{f^{(1)}}(\xi_1)\overline{(-i\xi_1)^\alpha\widehat{f^{(1)}}(\xi_1)}d\xi_1\nonumber&\\
&+\int^{+\infty}_0 (i\xi_1)^\alpha\widehat{f^{(1)}}(\xi_1)\overline{(-i\xi_1)^\alpha\widehat{f^{(1)}}(\xi_1)}d\xi_1\nonumber&\\
=&\cos(\pi\alpha)\int_{-\infty}^{+\infty}(i\xi_1)^\alpha\widehat{f^{(1)}}(\xi_1) \overline{(i\xi_1)^\alpha\widehat{f^{(1)}}(\xi_1)}d\xi_1\nonumber&\\
&+i\sin(\pi\alpha)\bigg(\int_0^{+\infty}(i\xi_1)^\alpha\widehat{f^{(1)}}(\xi_1) \overline{(i\xi_1)^\alpha\widehat{f^{(1)}}(\xi_1)}d\xi_1\nonumber&\\
&-\int_{-\infty}^{0}(i\xi_1)^\alpha\widehat{f^{(1)}}(\xi_1) \overline{(i\xi_1)^\alpha\widehat{f^{(1)}}(\xi_1)}d\xi_1\bigg),\nonumber&
\end{align}
where  $\widehat{ f^{(1)}}(\xi_1)=\frac{1}{(2\pi)^{\frac{1}{2}}}\int_{-\infty}^{+\infty} f^{(1)}(x,t)e^{-it\xi_1}dt$.

What's more, by $\overline{\widehat{f}(-\xi_1)}=\widehat{f}(\xi_1)$, we obtain that
\begin{align}
\int_0^{+\infty}(i\xi_1)^\alpha\widehat{f^{(1)}}(\xi_1) \overline{(i\xi_1)^\alpha\widehat{f^{(1)}}(\xi_1)}d\xi_1=\int_{-\infty}^{0}(i\xi_1)^\alpha\widehat{f^{(1)}}(\xi_1) \overline{(i\xi_1)^\alpha\widehat{f^{(1)}}(\xi_1)}d\xi_1.\nonumber
\end{align}

Therefore, we have
\begin{align}\label{title}
\int_{-\infty}^{+\infty}\widehat{(D^\alpha_1 f)}(\xi_1) \overline{\widehat{(D^{\alpha *}_1f)}(\xi_1)}d\xi_1
=\cos(\pi\alpha)\int_{-\infty}^{+\infty}(i\xi_1)^\alpha\widehat{f^{(1)}}(\xi_1) \overline{(i\xi_1)^\alpha\widehat{f^{(1)}}(\xi_1)}d\xi_1.
\end{align}

Substitute (\ref{title}) into (\ref{hjkljjnn}) yields
\begin{align}\label{hjkljjnnkm}
\int_{\mathbb{R}^3} D^\alpha_1 f(x) D^{\alpha *}_1f(x)dx
=&\cos(\pi\alpha)\int_{-\infty}^{+\infty}\int_{-\infty}^{+\infty}\left[\int_{-\infty}^{+\infty}(i\xi_1)^\alpha\widehat{f^{(1)}}(\xi_1) \overline{(i\xi_1)^\alpha\widehat{f^{(1)}}(\xi_1)}d\xi_1\right]dx_2dx_3\nonumber&\\
=&\cos(\pi\alpha)\int_{-\infty}^{+\infty}\int_{-\infty}^{+\infty}\left[\int_{-\infty}^{+\infty}\left|(i\xi_1)^\alpha\widehat{f^{(1)}}(\xi_1)\right|^2 d\xi_1\right]dx_2dx_3\nonumber&\\
=&\cos(\pi\alpha)\int_{-\infty}^{+\infty}\int_{-\infty}^{+\infty}\left[\int_{-\infty}^{+\infty}\left|D^\alpha_1 f(x)\right|^2 dx_1\right]dx_2dx_3\nonumber&\\
=&\cos(\pi\alpha)\|D^\alpha_1 f\|^2_{L^2(\mathbb{R}^3)}.&
\end{align}

In a similar way, we can prove $\int_{\mathbb{R}^3} D^\alpha_j f(x) D^{\alpha *}_jf(x)dx=\cos(\pi\alpha)\|D^\alpha_j f\|^2_{L^2(\mathbb{R}^3)}(j=2,3)$, which concludes (\ref{eq-newlem}).

\begin{lemma}\label{REq}
Assume $\alpha>0$ and $\alpha\neq m+\frac{1}{2}$, $m\in \mathbb{N}$, then $F_{L}^\alpha(\mathbb{R}^3)$, $F_{R}^\alpha(\mathbb{R}^3)$, $F_{S}^\alpha(\mathbb{R}^3)$ and $H^\alpha(\mathbb{R}^3)$ are equivalent.
\end{lemma}
{\bf Proof.} By Lemma \ref{newlem}, $F_{S}^\alpha(\mathbb{R}^3)$ and $F_{L}^\alpha(\mathbb{R}^3)$ are equivalent. On the other hand, by Lemma \ref{lem2-3}, we prove this Lemma.

\begin{lemma}\label{REq1}
Assume $\alpha>0$ and $\alpha\neq m+\frac{1}{2}$, $m\in\mathbb{N}$, then $F_{S,0}^\alpha(\Omega)$ and $H^\alpha_0(\Omega)$ are equivalent.
\end{lemma}

{\bf Proof.} Let $\tilde{f}$ be extension of $f\in C_0^\infty(\Omega)$ by zero outside $\Omega$, then $\mathrm{supp}(\tilde{f})\subset \Omega$, where $\mathrm{supp}(\tilde{f})=\{x: \tilde{f}(x)\neq 0\}$. What's more, $\mathrm{supp}(D_1^\alpha\tilde{f})\subset (a_1,+\infty)\times(a_2,b_2)\times(a_3,b_3)$, $\mathrm{supp}(D_2^\alpha\tilde{f})\subset (a_1,b_1)\times(a_2,+\infty)\times(a_3,b_3)$, $\mathrm{supp}(D_3^\alpha\tilde{f})\subset (a_1,b_1)\times(a_2,b_2)\times(a_3,+\infty)$,
 $\mathrm{supp}(D_1^{\alpha*}\tilde{f})\subset (-\infty,b_1)\times(a_2,b_2)\times(a_3,b_3)$, $\mathrm{supp}(D_2^{\alpha*}\tilde{f})\subset (a_1,b_1)\times(-\infty,b_2)\times(a_3,b_3)$, $\mathrm{supp}(D_3^{\alpha*}\tilde{f})\subset (a_1,b_1)\times(a_2,b_2)\times(-\infty,b_3)$.

Therefore, $\mathrm{supp}(D^\alpha \tilde{f}\cdot D^{\alpha *} \tilde{f})\subset \Omega$. This implies,
\begin{align}
 \|f\|_{F_{S,0}^\alpha(\Omega)}=\|\tilde{f}\|_{F_{S,0}^\alpha(\mathbb{R}^3)}, \|f\|_{H^\alpha_0(\Omega)}=\|\tilde{f}\|_{H^\alpha(\mathbb{R}^3)}.
 \end{align}

By Lemma \ref{REq}, we obtain that $F_{S,0}^\alpha(\Omega)$ and $H^\alpha_0(\Omega)$ are equivalent.\qed

Based on above Lemmas, we give the main result of this section.

\begin{thm}\label{the2-1}
Assume $\alpha>0$ and $\alpha\neq m+\frac{1}{2}$, $m\in\mathbb{N}$, then $F_{L,0}^\alpha(\Omega)$, $F_{R,0}^\alpha(\Omega)$ and  $H^\alpha_0(\Omega)$ are equivalent.
\end{thm}
{\bf Proof.} Let $\tilde{f}$ be extension of $f\in C_0^\infty(\Omega)$ by zero outside $\Omega$, then
\begin{align}
\|f\|_{F_{L,0}^\alpha(\Omega)}^2=&\|f\|_{L^2(\Omega)}^2+\|\nabla^\alpha f\|_{L^2(\Omega)}^2
\leq\|\tilde{f}\|_{L^2(\mathbb{R}^3)}^2+\|D^\alpha \tilde{f}\|_{L^2(\mathbb{R}^3)}^2\nonumber&\\
=&\|\tilde{f}\|_{F_{L}^\alpha(\mathbb{R}^3)}^2
\leq C\|\tilde{f}\|_{H^\alpha(\mathbb{R}^3)}^2=C\|f\|_{H^\alpha_0(\Omega)}^2.\nonumber&
\end{align}
That is, $H^\alpha_0(\Omega)\subseteq F_{L,0}^\alpha(\Omega)$.

On the other hand, by Lemma \ref{REq}, Lemma \ref{REq1} and Cauchy inequality\cite{E} $ab\leq \frac{a^2}{4\varepsilon}+\varepsilon b^2$ $(\varepsilon>0)$,
\begin{align}\label{inequality}
|f|_{H^\alpha(\Omega)}^2\leq& C|f|_{F_{S,0}^\alpha(\Omega)}^2=C\left|\int_\Omega \nabla^\alpha f(x)\cdot \nabla^{\alpha*}f(x)dx\right|\nonumber&\\
\leq&\frac{C}{4\varepsilon}\|\nabla^\alpha f\|_{L^2(\Omega)}^2+C\varepsilon\|\nabla^{\alpha *} f\|_{L^2(\Omega)}^2\nonumber&\\
=&\frac{C}{4\varepsilon}|f|_{F_{L,0}^\alpha(\Omega)}^2+C\varepsilon|f|_{F_{R,0}^\alpha(\Omega)}^2\nonumber
\leq\frac{C}{4\varepsilon}|f|_{F_{L,0}^\alpha(\Omega)}^2+C\varepsilon|\tilde{f}|_{F_{R,0}^\alpha(\mathbb{R}^3)}^2\nonumber&\\
\leq&\frac{C}{4\varepsilon}|f|_{F_{L,0}^\alpha(\Omega)}^2+C_1\varepsilon|\tilde{f}|_{H^\alpha_0(\mathbb{R}^3)}^2
=\frac{C}{4\varepsilon}|f|_{F_{L,0}^\alpha(\Omega)}^2+C_1\varepsilon|f|_{H^\alpha_0(\Omega)}^2.&
\end{align}

Let $\varepsilon=\frac{1}{2C_1}$ in (\ref{inequality}), then
\begin{align}
|f|_{H^\alpha(\Omega)}^2\leq CC_1|f|_{F_{L,0}^\alpha(\Omega)}^2.
\end{align}

That is, $F_{L,0}^\alpha(\Omega)\subseteq H^\alpha_0(\Omega)$.

Now, we conclude  $F_{L,0}^\alpha(\Omega)=H^\alpha_0(\Omega)$.

Similarly, we can prove that $F_{R,0}^\alpha(\Omega)$ and  $H^\alpha_0(\Omega)$ are equivalent.

\section{Existence of solution to (\ref{eq1-4})}

\begin{defn}\label{defi3-1}
(cf. \cite{ER})For $\alpha>0$, $x\in[a,b]$ and function $f:[a,b]\rightarrow \mathbb{R}$, define
\begin{align}
D^{-\alpha}_{[a,x]}f(x)=\frac{1}{\Gamma(\alpha)}\int_a^x(x-t)^{\alpha-1}f(t)dt, D^{-\alpha*}_{[x,b]}f(x)=\frac{1}{\Gamma(\alpha)}\int_x^b(t-x)^{\alpha-1}f(t)dt.\nonumber
\end{align}
\end{defn}

By Appendix in \cite{ER}, we know that
\begin{align}\label{eq3-1}
D^{-\alpha}_{[a,x]}D^{\alpha}_{[a,x]} f(x)=f(x),\ \  D^{-\alpha*}_{[x,b]}D^{\alpha*}_{[x,b]} f(x)=f(x).
\end{align}

Based on Definition \ref{defi3-1}, for $x=(x_1,x_2,x_3)\in\Omega$, $i=1,2,3$ and function $f:\Omega\rightarrow \mathbb{R}$, define
\begin{align}
D^{-\alpha}_{[a_i,x_i]}f(x)=\frac{1}{\Gamma(\alpha)}\int_{a_i}^{x_i}(x_i-t)^{\alpha-1}f^{(i)}(x,t)dt, D^{-\alpha*}_{[x_i,b_i]}f(x)=\frac{1}{\Gamma(\alpha)}\int_{x_i}^{b_i}(t-x_i)^{\alpha-1}f^{(i)}(x,t)dt.\nonumber
\end{align}

By (\ref{eq3-1}), we know that
\begin{align}\label{eq3-3}
D^{-\alpha}_{[a_i,x_i]}D^{\alpha}_{[a_i,x_i]} f(x)=f(x),\ \  D^{-\alpha*}_{[x_i,b_i]}D^{\alpha*}_{[x_i,b_i]} f(x)=f(x)(i=1,2,3).
\end{align}

Based on above definitions and notations, we obtain the following property of operators $D^{-\alpha}_{[a_i,x_i]}$ and $D^{-\alpha*}_{[x_i,b_i]}(i=1,2,3)$.

\begin{lemma}\label{lemm3-1}
For $\alpha>0$, $x=(x_1,x_2,x_3)\in \Omega$ and function  $f\in L^2(\Omega)$, there exists a constant $C=C(\alpha,\Omega)$ such that
\begin{align}
\|D^{-\alpha}_{[a_i,x_i]}f\|_{L^2(\Omega)}\leq C\|f\|_{L^2(\Omega)}, \|D^{-\alpha*}_{[x_i,b_i]}f\|_{L^2(\Omega)}\leq C\|f\|_{L^2(\Omega)}\ \ (i=1,2,3)\nonumber.
\end{align}
\end{lemma}
{\bf Proof.} Since $D_{[a_1,x_1]}^{-\alpha}f(x)=\frac{x_1^{\alpha-1}}{\Gamma(\alpha)}*f(x_1,x_2,x_3)$, by Young's inequality\cite{E}
\begin{align}\label{eq33}
\|v*w\|_{L^2([a_1,b_1])}^2\leq\|v\|_{L^1([a_1,b_1])}^2\|w\|_{L^2([a_1,b_1])}^2,
\end{align}
it yields,
\begin{align}
\|D_{[a_1,x_1]}^{-\alpha}f\|^2_{L^2(\Omega)}=&\int_{a_1}^{b_1}\int_{a_2}^{b_2}\int_{a_3}^{b_3}\left|D_{[a_1,x_1]}^{-\alpha}f(x_1,x_2,x_3)\right|^2dx_3dx_2dx_1\nonumber&\\
=&\int_{a_2}^{b_2}\int_{a_3}^{b_3}\left[\int_{a_1}^{b_1}\left|D_{[a_1,x_1]}^{-\alpha}f(x_1,x_2,x_3)\right|^2dx_1\right]dx_3dx_2\nonumber&\\
\leq&\int_{a_2}^{b_2}\int_{a_3}^{b_3}\left[\frac{1}{(\Gamma(\alpha))^2}\int_{a_1}^{b_1}|f(x_1,x_2,x_3)|^2dx_1\left(\int_{a_1}^{b_1}|x_1^{\alpha-1}|dx_1\right)^2\right]dx_3dx_2\nonumber&\\
\leq&\int_{a_2}^{b_2}\int_{a_3}^{b_3}\left[\frac{(|a_1|^\alpha+|b_1|^\alpha)^2}{(\Gamma(\alpha+1))^2}\int_{a_1}^{b_1}|f(x_1,x_2,x_3)|^2dx_1\right]dx_3dx_2\nonumber&\\
=&\frac{(|a_1|^\alpha+|b_1|^\alpha)^2}{(\Gamma(\alpha+1))^2}\int_{a_1}^{b_1}\int_{a_2}^{b_2}\int_{a_3}^{b_3}|f(x_1,x_2,x_3)|^2dx_3dx_2dx_1\nonumber&\\
=&\frac{(|a_1|^\alpha+|b_1|^\alpha)^2}{(\Gamma(\alpha+1))^2}\|f\|^2_{L^2(\Omega)}.&
\end{align}

Similarly, $\|D_{[a_i, x_i]}^{-\alpha}f\|^2_{L^2(\Omega)}\leq\frac{(|a_i|^\alpha+|b_i|^\alpha)^2}{(\Gamma(\alpha+1))^2}\|f\|^2_{L^2(\Omega)}\ (i=2,3)$.

Let $C=\max\limits_{i=1,2,3}\left\{\frac{|a_i|^\alpha+|b_i|^\alpha}{\Gamma(\alpha+1)}\right\}$, then
\begin{align}
\|D^{-\alpha}_{[a_i, x_i]}f\|_{L^2(\Omega)}\leq C\|f\|_{L^2(\Omega)}\ \ (i=1,2,3).
\end{align}

In a similar way, we can prove that
\begin{align}
\|D^{-\alpha*}_{[x_i,b_i]}f\|_{L^2(\Omega)}\leq C\|f\|_{L^2(\Omega)}\ \ (i=1,2,3).
\end{align}

\begin{lemma}\label{lemm3-2}
Assume $\alpha>0$, $u\in F^\alpha_{L,0}(\Omega)$ and $v\in F^\alpha_{R,0}(\Omega)$, then there exists a constant $C=C(\alpha,\Omega)$ such that
\begin{align}\label{eq3-4}
\|u\|_{L^2(\Omega)}\leq C|u|_{F^\alpha_{L,0}(\Omega)},\ \ \|v\|_{L^2(\Omega)}\leq C|v|_{F^\alpha_{R,0}(\Omega)}.
\end{align}

\end{lemma}
{\bf Proof.} It follows from Lemma \ref{lemm3-1} and (\ref{eq3-3}) that
\begin{align}\label{hja}
\|u\|_{L^2(\Omega)}=\|D^{-\alpha}_{[a_j,x_j]}D^{\alpha}_{[a_j,x_j]} u\|_{L^2(\Omega)}\leq \bar{C}\left\|\frac{\partial^\alpha u}{\partial x_j^\alpha}\right\|_{L^2(\Omega)}(j=1,2,3),
\end{align}
where $\bar{C}=\bar{C}(\alpha,\Omega)$.

It follows from (\ref{hja}) that
\begin{align}
\|u\|_{L^2(\Omega)}^2\leq \frac{\bar{C}^2}{3}\sum_{j=1}^3\left\|\frac{\partial^\alpha u}{\partial x_j^\alpha}\right\|_{L^2(\Omega)}^2=\frac{\bar{C}^2}{3}|u|^2_{F^\alpha_{L,0}(\Omega)}.
\end{align}

Let $C=\frac{\bar{C}}{\sqrt{3}}$, this concludes the first equation of (\ref{eq3-4}).

Similarly, we can prove the second equation of (\ref{eq3-4}).

\begin{lemma}\label{lemm3-3}
Assume $\alpha>2.5$, $\alpha\neq m+0.5$, $m\in \mathbb{N}$ and $u\in F^\alpha_{L,0}(\Omega)$, then there exists a constant $K=K(\alpha,\Omega)$ such that
\vskip1mm {\bf (i).} $\|u(x)-u(y)\|\leq K\|\nabla^\alpha u\|_{L^2(\Omega)}\|x-y\|$.
\vskip1mm {\bf (ii).} $\|\nabla u(x)-\nabla u(y)\|\leq K\|\nabla^\alpha u\|_{L^2(\Omega)}\|x-y\|^\lambda$, where $0<\lambda\leq\alpha-1.5-[\alpha-1.5]$ as $2.5<\alpha<3.5$ and $\lambda=1$ as $\alpha\geq3.5$.
\end{lemma}
{\bf Proof.}{\bf (i).} By Theorem \ref{the2-1}, $u\in F^\alpha_{L,0}(\Omega)=H^\alpha_{0}(\Omega)\hookrightarrow C^{1}(\Omega)$(cf. [2, Theorem 4.57]).  This implies
\begin{align}\label{eq3-5}
{\|u(x)-u(y)\|}=&\|\nabla u(\xi)\cdot{(x-y)}\|\leq \|u\|_{C^1(\Omega)}\|x-y\|\leq C_1\|u\|_{H^\alpha_0(\Omega)}\|x-y\|\nonumber&\\
\leq& C_2\|u\|_{F^\alpha_{L,0}(\Omega)}\|x-y\|.&
\end{align}

On the other hand, by Lemma \ref{lemm3-2}, we know there exists a constant $C_3=C_3(\alpha,\Omega)$ such that
\begin{align}\label{eq3-6}
\|u\|_{F^\alpha_{L,0}(\Omega)}\leq C_3\|\nabla^\alpha u\|_{L^2(\Omega)}.
\end{align}

By (\ref{eq3-5}) and (\ref{eq3-6}), we obtain that
\begin{align}\label{eq3-7}
{\|u(x)-u(y)\|}\leq K\|\nabla^\alpha u\|_{L^2(\Omega)}{\|x-y\|}.
\end{align}

This concludes {\bf (i)}.

\vskip1mm {\bf (ii).} Here we divide our discussion into two different cases:

 {\bf Case 1($2.5<\alpha<3.5$).}
 Since $F^\alpha_{L,0}(\Omega)=H^\alpha_{0}(\Omega)\hookrightarrow C^{1,\lambda}(\Omega)$ $(0<\lambda<\alpha-1.5-[\alpha-1.5])$[2, Theorem 4.57], then
\begin{align}\label{eq3-8}
\frac{\|\nabla u(x)-\nabla u(y)\|}{\|x-y\|^\lambda}\leq \|u\|_{C^{1,\lambda}(\Omega)}\leq C\|u\|_{F^\alpha_{L,0}(\Omega)}\leq K\|\nabla^\alpha u\|_{L^2(\Omega)}.
\end{align}
\vskip1mm {\bf Case 2($\alpha\geq3.5$).} Since $F^\alpha_{L,0}(\Omega)=H^\alpha_{0}(\Omega)\hookrightarrow C^{2}(\Omega)$[2, Theorem 4.57], then
\begin{align}\label{eq3-9}
\frac{\|\nabla u(x)-\nabla u(y)\|}{\|x-y\|}\leq \|u\|_{C^{2}(\Omega)}\leq C\|u\|_{F^\alpha_{L,0}(\Omega)}\leq K\|\nabla^\alpha u\|_{L^2(\Omega)}.
\end{align}

This concludes {\bf (ii)}.

By (\ref{eq1-4}), we know $H$ is a functional about $v$ and $\eta$, where $v$ and $\eta$ are constrained by (\ref{eq1-3}). In this paper, we write $H$ as a functional only about $v$. Therefore, (\ref{eq1-3}) should admit a unique solution. Otherwise, the definition of functional $H$ is ambiguous. As to the well-define of $H$, we have the following result.
\begin{lemma}\label{lemm3-4}
Assume $v(x,s)\in\mathcal{F}$ and $v(\cdot,s)|_{\mathbb{R}^3\setminus\Omega}=0$ for each $s\in[0,\tau]$, then for each $x\in\bar{\Omega}$, (\ref{eq1-3}) admits a unique solution $\eta(s;t,x)\in C([0,\tau],\bar{\Omega})$. Moreover, for each $s\in [0,\tau]$, $\eta(s;t,x)\in [C^{[\alpha-1.5],\lambda}(\Omega)]^3$, where $0<\lambda\leq\alpha-1.5-[\alpha-1.5]$.
\end{lemma}
{\bf Proof.} Based on {\bf (i)} in Lemma \ref{lemm3-3}, this conclusion can be proved in a similar way with Lemma 2.2 in \cite{HH}.

As to the existence of $h^{-1}:\Omega\rightarrow\Omega$, we have the following result.
\begin{lemma}\label{lemm3-5}
Assume $v(x,s)\in\mathcal{F}$ and $v(\cdot,s)|_{\mathbb{R}^3\setminus\Omega}=0$ for each $s\in[0,\tau]$, then for $s,t\in[0,\tau]$ and $x\in\bar{\Omega}$, $\eta(s;t,x)$ defined by (\ref{eq1-3}) is differential with respect to $x$. Define $\Theta(s;t,x)\triangleq\nabla_x \eta(s;t,x)$, then \\
{\bf (i).}$\Theta(s;t,x)$ is the solution of
 \begin{equation}\label{fenshueq-hj5}\left\{\begin{aligned}
 &\frac{d\Theta(s;t,x)}{ds}=\nabla_\eta v(\eta(s;t,x),s)\Theta (s;t,x), &\\
 &\Theta (t;t,x)=I, &\\
 \end{aligned}\right.\end{equation}
and
\begin{align}
\det(\Theta(s;t,x))=e^{\int_t^s\sum\limits_{i=1}^3 v_{i,\eta_i}(\eta(r;t,x),r)dr},
\end{align}
where $v_{i,\eta_i}(\eta(r;t,x),r)=\frac{\partial v_i(\eta(r;t,x),r)}{\partial \eta_i}(i=1,2,3)$.\\
{\bf (ii).} Define $h(x)\triangleq\eta(0;\tau,x)$ as (\ref{eq1-3}), then $h$ is a $1$-to-$1$ and onto mapping which ensures the existence of $h^{-1}(x)$ in (\ref{eq1-2real}).
\end{lemma}
{\bf Proof.} {\bf (i).}Based on Theorem \ref{the2-1}, this conclusion can be obtained in a similar way with Lemma 2.3 in \cite{HH}.\\
{\bf (ii).} For any $x\in\Omega$, $\det{(\nabla_x{{h}}(x))}=e^{-\int_0^\tau\sum\limits_{i=1}^3v_{i,\eta_i}(\eta(s;\tau,x),s)ds}\neq 0$, by Inverse Function Theorem[5, Theorem 7 in Appendix C], we know ${{h}}$ is a $1$-to-$1$ and onto mapping. This ensures the existence of $h^{-1}(x)$ in (\ref{eq1-2real}).

\begin{lemma}\label{lemm3-6}
Assume $\{v_n(x,s)\}$ is a bounded sequence on $\mathcal{F}$ with
\begin{align}\label{{eq3-10}}
\|v_n\|_{\mathcal{F}}^2=\int_0^\tau\|\nabla^\alpha v_n(\cdot,s)\|_{L^2(\Omega)}^2ds\leq M<+\infty,
\end{align}
and $v_n(\cdot,s)|_{\mathbb{R}^3\setminus\Omega}=0$ for each $s\in[0,\tau]$, then
\vskip1mm {\bf(i).}$\{v_n(x,s)\}$ is a weakly compact set on $\mathcal{F}$.
\vskip1mm {\bf(ii).}If we denote $n_k$ as the sequence number of a weakly convergent subsequence $\{v_{n_k}(x,s)\}$ with weak limit $v(x,s)$, then
\begin{align}\label{{eq3-11}}
M \geq\lim_{n_k\rightarrow \infty}\inf \int_0^\tau\|\nabla^\alpha v_{n_k}(\cdot,s)\|_{L^2(\Omega)}^2ds\geq \int_0^\tau\|\nabla^\alpha v(\cdot,s)\|_{L^2(\Omega)}^2ds.
\end{align}
\vskip1mm {\bf(iii).} Let $v_{n_k}(x,s)$, $v(x,s)$ be functions defined in {\bf(ii)}. Consider the equations
\begin{align}\label{eq3-d12}
 \frac{{d\eta}_{n_k}(s;t,x)}{ds}=v_{n_k}(\eta_{n_k}(s;t,x),s) \  \ \mathrm{with}\ \ \eta_{n_k}(t;t,x)=x,
 \end{align}
 and
 \begin{align}\label{eq3-13}
 \frac{d\eta(s;t,x)}{ds}=v(\eta(s;t,x),s) \  \ \mathrm{with}\ \ \eta(t;t,x)=x,
 \end{align}
 then for each $(x,t)\in\Omega\times[0,\tau]$, these two equations have a unique solution $\eta_{n_k}(s;t,x)$,$\eta(s;t,x)\in C([0,\tau], \bar{\Omega})$, respectively. Furthermore, for each $s\in [0,\tau]$, $\eta_{n_k}(s;t,x),\eta(s;t,x)\in [C^{[\alpha-1.5],\lambda}(\Omega)]^3$ with $\eta_{n_k}(s;t,x)\xrightarrow{k} {\eta}(s;t,x)$ uniformly on $[0,\tau]$, where $0<\lambda\leq\alpha-1.5-[\alpha-1.5]$.
\end{lemma}

{\bf Proof.}{\bf(i).} By Lemma \ref{lemm3-2}, we know $\|\cdot\|_\mathcal{F}$ is a norm and $\mathcal{F}$ is a separable Hilbert space. This implies {\bf(i)} for the fact that any closed ball in a separable Hilbert space is a weakly compact set.

{\bf(ii).} Since $\|\cdot\|_\mathcal{F}$ is a norm, by the lower weak semi-continuity of norm, we obtain (\ref{{eq3-11}}).

{\bf(iii).} By Lemma \ref{lemm3-4}, we know that equations (\ref{eq3-d12}) and (\ref{eq3-13}) have a unique solution $\eta_{n_k}(s;t,x)$,$\eta(s;t,x)\in [C^{[\alpha-1.5],\lambda}(\Omega)]^3$ respectively, with $\eta_{n_k}(s;t,x)\xrightarrow{k} {\eta}(s;t,x)$ uniformly on $[0,\tau]$, where $0<\lambda\leq\alpha-1.5-[\alpha-1.5]$. \qed

Before we give a proof of Theorem \ref{th-A10}, let's recall the following Lemma.
\begin{lemma}\label{mainlemm}
(cf. [9, Theorem 1.A.4])Let $E$ be a weakly compact set on Banach space $X$. If $H:E\rightarrow \mathbb{R}$ is a lower weak semi-continuous$(\mathrm{l.w.c\ in short)}$ functional, then there exists $v_0\in E$ such that $H(v_0)=\inf\limits_{v\in E} H(v)$.
\end{lemma}

{\bf Proof of Theorem \ref{th-A10}.} Define function $\tilde{v}(x,s)\equiv 0$ on $\Omega\times [0,\tau]$. Since $\tilde{v}(x,s)\in\mathcal{F}$ and $H(\tilde{v})=\|T(\cdot)-D(\cdot)\|_{L^2(\Omega)}^2\leq G|\Omega|\triangleq M$, we only need to show the existence of global minimizer of $H(v)$ on the ball
\begin{align}
B_M\triangleq\{v(x,s): \|v\|_\mathcal{F}^2\leq M\}.
\end{align}

By {\bf (i)} in Lemma \ref{lemm3-6}, we know $B_M$ is a weakly compact set. Choose $\{v_n\}\in B_M$, then there exists a weakly convergent subsequence $\{v_{n_k}\}$ such that
\begin{align}
v_{n_k}\xrightharpoonup{k} v \in B_M.
\end{align}

By {\bf (ii)} in Lemma \ref{lemm3-6}, we know that
\begin{align}\label{{neweq3-11}}
M \geq\lim_{n_k\rightarrow \infty}\inf \int_0^\tau\|\nabla^\alpha v_{n_k}(\cdot,t)\|_{L^2(\Omega)}^2dt\geq \int_0^\tau\|\nabla^\alpha v(\cdot,t)\|_{L^2(\Omega)}^2dt.
\end{align}

By {\bf (iii)} in Lemma \ref{lemm3-6}, $\eta_{n_k}(s;t,x)\xrightarrow{k} {\eta}(s;t,x)$ for all $(x,s,t)\in\Omega\times[0,\tau]\times[0,\tau]$,  where $\eta_{n_k}(s;t,x)$, ${\eta}(s;t,x)$ are the solution of (\ref{eq3-d12}) and (\ref{eq3-13}), respectively. Furthermore, we have
\begin{align}\label{eq3-24}
h_{n_k}(x)=\eta_{n_k}(0;\tau,x)\xrightarrow{k}\eta(0;\tau,x)=h(x) \ \mathrm{for\ all}\ x\in\Omega.
\end{align}

Define $\Theta_{n_k}(s;0,x)\triangleq\nabla_x\eta_{n_k}(s;0,x)$, $\Theta(s;0,x)\triangleq\nabla_x\eta(s;0,x)$, then by Lemma \ref{lemm3-5}, $\Theta_{n_k}(s;0,x)$ and $\Theta(s;0,x)$ are the solutions of
\begin{equation}\label{newfenshueq-hj5}\left\{\begin{aligned}
 &\frac{d\Theta_{n_k}(s;0,x)}{ds}=\nabla_{\eta_{n_k}} v_{n_k}(\eta_{n_k}(s;0,x),s)\Theta_{n_k} (s;0,x), &\\
 &\Theta_{n_k} (0;0,x)=I, &\\
 \end{aligned}\right.\end{equation}
and
\begin{equation}\label{new1fenshueq-hj5}\left\{\begin{aligned}
 &\frac{d\Theta(s;0,x)}{ds}=\nabla_\eta v(\eta(s;0,x),s)\Theta (s;0,x), &\\
 &\Theta(0;0,x)=I, &\\
 \end{aligned}\right.\end{equation}
respectively.

Based on these notations, here we claim that the functional
\begin{align}\label{claim}
 H:B_M\rightarrow \mathbb{R}\ \ \mathrm{is\ l.w.c}.
 \end{align}

  The proof of claim (\ref{claim}) can be divided into following five steps.

{\bf Step 1.} We claim that there exists a constant $0<\tilde{M}<+\infty$ such that
\begin{align}
\|\Theta_{n_k}(s;0,x)\|\leq \tilde{M},\ \ \|\Theta(s;0,x)\|\leq \tilde{M}.
\end{align}

By (\ref{newfenshueq-hj5}), we have
\begin{align}
\Theta_{n_k}(s;0,x)=I+\int_0^s \nabla_{\eta_{n_k}} v_{n_k}(\eta_{n_k}(r;0,x),r)\Theta_{n_k} (r;0,x)dr,
\end{align}
this yields
\begin{align}
\|\Theta_{n_k}(s;0,x)\|=\|I\|+\int_0^s \|\nabla_{\eta_{n_k}} v_{n_k}(\eta_{n_k}(r;0,x),r)\|\|\Theta_{n_k} (r;0,x)\|dr.\nonumber
\end{align}

By Grownwall inequality and Lemma \ref{lemm3-2}, we have
\begin{align}
\|\Theta_{n_k}(s;0,x)\|\leq&\|I\|e^{\int_0^s\|\nabla_{\eta_{n_k}} v_{n_k}(\eta_{n_k}(r;0,x),r)\|dr}
\leq\|I\|e^{\int_0^s\|v_{n_k}(\cdot,r)\|_{C^1(\Omega)}dr}\nonumber&\\
\leq&\|I\|e^{\int_0^sC\|v_{n_k}(\cdot,r)\|_{F^\alpha_{L,0}(\Omega)}dr}
\leq\|I\|e^{\int_0^s\tilde{C}\|\nabla^\alpha v_{n_k}(\cdot,r)\|_{L^2(\Omega)}dr}\nonumber&\\
\leq&\|I\|e^{\tilde{C}\tau^{\frac{1}{2}}(\int_0^s\|\nabla^\alpha v_{n_k}(\cdot,r)\|^2_{L^2(\Omega)}dr)^{\frac{1}{2}}}
\leq\|I\|e^{\tilde{C}{(\tau M)}^{\frac{1}{2}}}\triangleq \tilde{M}.\nonumber&
\end{align}

Similarly, we can prove that $\|\Theta(s;0,x)\|\leq \tilde{M}$.

{\bf Step 2.} We claim that $\nabla_x\left[\int_0^s [v_{n_k}(x,r)-v(x,r)]dr\right]\xrightarrow{k} O_{3\times 3}\ \ \mathrm{uniformly \ on\ } \Omega\times[0,\tau]$, where $O_{3\times 3}$ is a ${3\times 3}$ matrix whose elements are all zero.

Let $w_{n_k}(x,s)\triangleq v_{n_k}(x,s)-v(x,s)$, $z_{n_k}(x,s)\triangleq\int_0^sw_{n_k}(x,r)dr$, then $w_{n_k}\xrightharpoonup{k}0$ with
$\|w_{n_k}\|_\mathcal{F}\leq 2M$.

If $x,y\in\Omega$ and $s,t\in[0,\tau]$, then by {\bf (i)} in Lemma \ref{lemm3-3},

\begin{align}
 &\ \|z_{n_k}(x,s)-z_{n_k}(y,t)\|\nonumber&\\
 =&\left\|\int_0^s [w_{n_k}(x,r)-w_{n_k}(y,r)]dr\right\|+\left\|\int_s^t w_{n_k}(y,r)dr\right\| \nonumber &\\
 \leq &K\int_0^\tau \|\nabla^\alpha v_{n_k}(\cdot,r)\|_{L^2(\Omega)}dr\|x-y\|\nonumber&\\
 &+K\int_0^\tau \|\nabla^\alpha v(\cdot,r)\|_{L^2(\Omega)}dr\|x-y\|+\left\|\int_s^t w_{n_k}(y,r)dr\right\| \nonumber &\\
 \leq &K\tau^{\frac{1}{2}}\left(\int_0^\tau \|\nabla^\alpha v_{n_k}(\cdot,r)\|_{L^2(\Omega)}^2dr\right)^{\frac{1}{2}}\|x-y\|\nonumber &\\
 &+K\tau^{\frac{1}{2}}\left(\int_0^\tau \|\nabla^\alpha v(\cdot,r)\|_{L^2(\Omega)}^2dr\right)^{\frac{1}{2}}\|x-y\| +\left\|\int_s^t w_{n_k}(y,r)dr\right\| \nonumber & \\
 \leq &2KM^{\frac{1}{2}}\tau^{\frac{1}{2}}\|x-y\|+\left\|\int_s^t w_{n_k}(y,r)dr\right\| \nonumber & \\
  \leq& 2KM^{\frac{1}{2}}\tau^{\frac{1}{2}}\|x-y\|+C|t-s|^{\frac{1}{2}}\|w_{n_k}\|_\mathcal{F} \nonumber & \\
  \leq& 2KM^{\frac{1}{2}}\tau^{\frac{1}{2}}\|x-y\|+2CM|t-s|^{\frac{1}{2}}. \nonumber &
\end{align}

By Arzela-Ascoli Theorem\cite{LZS}, $\{z_{n_k}\}$ is relative compact on $[C(\Omega\times[0,\tau])]^3$.

By \cite{DGM,HH}, $z_{n_k}\xrightarrow{k} 0$ uniformly on $\Omega\times [0,\tau]$.

On the other hand, by {\bf (ii)} in Lemma \ref{lemm3-3}, we have
\begin{align}
&\ \ \ \ \|\nabla_x z_{n_k}(x,s)-\nabla_y z_{n_k}(y,t)\|\nonumber&\\
 &=\left\|\int_0^s [\nabla_x w_{n_k}(x,r)-\nabla_y w_{n_k}(y,r)]dr\right\|+\left\|\int_s^t \nabla_y w_{n_k}(y,r)dr\right\| \nonumber &\\
 &\leq K\int_0^\tau \|\nabla^\alpha v_{n_k}(\cdot,r)\|_{L^2(\Omega)}dr\|x-y\|^{\lambda}\nonumber &\\
 &+K\int_0^\tau \|\nabla^\alpha v(\cdot,r)\|_{L^2(\Omega)}dr\|x-y\|^{\lambda}+\left\|\int_s^t \nabla_y w_{n_k}(y,r)dr\right\| \nonumber &\\
 &\leq K\tau^{\frac{1}{2}}\left(\int_0^\tau \|\nabla^\alpha v_{n_k}(\cdot,r)\|_{L^2(\Omega)}^2dr\right)^{\frac{1}{2}}\|x-y\|^{\lambda}\nonumber&
\end{align}

\begin{align}
 &+K\tau^{\frac{1}{2}}\left(\int_0^\tau \|\nabla^\alpha v(\cdot,r)\|_{L^2(\Omega)}^2dr\right)^{\frac{1}{2}}\|x-y\|^{\lambda}+\left\|\int_s^t \nabla_y w_{n_k}(y,r)dr\right\| \nonumber & \\
 &\leq 2KM^{\frac{1}{2}}\tau^{\frac{1}{2}}\|x-y\|^{\lambda}+\left\|\int_s^t \nabla_y w_{n_k}(y,r)dr\right\| \nonumber & \\
  &\leq 2KM^{\frac{1}{2}}\tau^{\frac{1}{2}}\|x-y\|^{\lambda}+|t-s|^{\frac{1}{2}}K\|w_{n_k}\|_{\mathcal{F}} \nonumber &\\
  &\leq 2KM^{\frac{1}{2}}\tau^{\frac{1}{2}}\|x-y\|^{\lambda}+2MK|t-s|^{\frac{1}{2}}. \nonumber &
\end{align}

Therefore, $\{\nabla_x z_{n_k}(x,s)\}$ is relative compact on $[C(\Omega\times[0,\tau])]^3$.

Now, choose any convergent subsequence of $\{\nabla_x z_{n_k}(x,s)\}$ with limit $P_{3\times 3}$. That is,
\begin{align}
\nabla_x z_{n_k}(x,s)\xrightarrow{k}P_{3\times 3}\  \mathrm{uniformly\ on\ } \Omega\times [0,\tau].
\end{align}

Since $z_{n_k}\xrightarrow{k} 0$ uniformly on $\Omega\times[0,\tau]$, then $P_{3\times 3}=\nabla_x 0=O_{3\times 3}$.

This implies, $\nabla_x\left[\int_0^s [v_{n_k}(x,r)-v(x,r)]dr\right]\xrightarrow{k} O_{3\times 3}\ \ \mathrm{uniformly \ on\ } \Omega\times[0,\tau]$.

{\bf Step 3.} We claim that $\|\Theta_{n_k}\Theta_{n_k}^T- \Theta\Theta^T\|\xrightarrow{k}0$. Here for the sake of simplicity,  we denote $\Theta_{n_k}(s;0,x)$ and $\Theta(s;0,x)$ by $\Theta_{n_k}$ and $\Theta$, respectively.

\begin{align}
&\|\Theta_{n_k}(s;0,x)-\Theta(s;0,x)\|\nonumber&\\
&\leq \int_0^s \|[\nabla_{\eta_{n_k}} v_{n_k}(\eta_{n_k}(r;0,x),r)-\nabla_{\eta} v_{n_k}(\eta(r;0,x),r)]\Theta_{n_k}(r;0,x)\|dr \nonumber&\\
&+\int_0^s\|\nabla_{\eta} v_{n_k}(\eta(r;0,x),r)[\Theta_{n_k}(r;0,x)-\Theta(r;0,x)]\|dr\nonumber& \\
&+\left\|\int_0^s  [\nabla_{\eta} v_{n_k}(\eta(r;0,x),r)-\nabla_{\eta} v(\eta(r;0,x),r)]\Theta(r;0,x)dr\right\| &\nonumber\\
&\leq \int_0^s K\|\nabla^\alpha v_{n_k}(\cdot,r)\|_{L^2(\Omega)}\|\eta_{n_k}(r;0,x)-\eta(r;0,x)\|^{\lambda}\|\Theta_{n_k}(r;0,x)\|dr \nonumber&\\
&+ \int_0^s\|\nabla_{\eta} v_{n_k}(\eta(r;0,x),r)\|\|\Theta_{n_k}(r;0,x)-\Theta(r;0,x)\| dr \nonumber& \\
&+ \left\|\int_0^s[\nabla_\eta v_{n_k}(\eta(r;0,x),r)-\nabla_\eta v(\eta(r;0,x),r)]\Theta(r;0,x) dr\right\| \nonumber&
\end{align}

\begin{align}\label{eq-HY}
&\leq K\tilde{M}\tau^{\frac{1}{2}}\left(\int_0^s \|\nabla^\alpha v_{n_k}(\cdot,r)\|_{L^2(\Omega)}^2dr\right)^{\frac{1}{2}}\|\eta_{n_k}(r;0,x)-\eta(r;0,x)\|^{\lambda}_{C([0,\tau]:\overline{\Omega})} \nonumber&\\
&+ \int_0^s\|\nabla_\eta v_{n_k}(\eta(r;0,x),r)\|\|\Theta_{n_k}(r;0,x)-\Theta(r;0,x)\| dr \nonumber& \\
&+ \left\|\int_0^s[\nabla_{\eta} v_{n_k}(\eta(r;0,x),r)-\nabla_\eta v(\eta(r;0,x),r)]\Theta(r;0,x) dr\right\|  \nonumber&\\
&\leq K\tilde{M}\tau^{\frac{1}{2}}M^{\frac{1}{2}}\|\eta_{n_k}(r;0,x)-\eta(r;0,x)\|^{\lambda}_{C([0,\tau]:\overline{\Omega})} \nonumber&\\
&+ \left\|\int_0^s[\nabla_\eta v_{n_k}(\eta(r;0,x),r)-\nabla_\eta v(\eta(r;0,x),r)]\Theta(r;0,x) dr\right\| \nonumber&\\
&+ \int_0^s\|\nabla_\eta v_{n_k}(\eta(r;0,x),r)\|\|\Theta_{n_k}(r;0,x)-\Theta(r;0,x)\| dr\nonumber&\\
&\triangleq I_1+I_2+I_3.&
\end{align}

On the other hand, by {\bf (iii)} in Lemma \ref{lemm3-6}, we have $\|\eta_{n_k}(s;0,x)-\eta(s;0,x)\|_{C([0,\tau]:\overline{\Omega})}\xrightarrow{k} 0$.

This leads to $I_1\xrightarrow{k} 0$ . By {\bf Step 2}, we know that $\nabla_x\left[\int_0^s v_{n_k}(x,r)-v(x,r)dr\right]\xrightarrow{k} O_{3\times 3}$ uniformly on $\Omega\times[0,\tau]$. This implies that $I_2\xrightarrow{k} 0$.
\vskip1mm Hence, $I_1+I_2\xrightarrow{k} 0$. Then there exist $N=N(\varepsilon)$, such that $I_1+I_2<\varepsilon$ as $n_k>N$. Therefore, it follows from (\ref{eq-HY})  that
\begin{align}
\|\Theta_{n_k}(s;0,x)-\Theta(s;0,x)\|&\leq \varepsilon+ \int_0^s\|\nabla_\eta v_{n_k}(\eta(r;0,x),r)\|\|\Theta_{n_k}(r;0,x)-\Theta(r;0,x)\| dr, \nonumber
\end{align}
and Gronwall inequality implies that
\begin{align}
\|\Theta_{n_k}(s;0,x)-\Theta(s;0,x)\|&\leq \varepsilon e^{\int_0^s\|\nabla_\eta v_{n_k}(\eta(r;0,x),r)\| dr} \leq \varepsilon e^{K\tau^{\frac{1}{2}}(\int_0^\tau \|\nabla^\alpha v_{n_k}(\cdot,r)\|^2_{L^2(\Omega)} dr)^{\frac{1}{2}}}\leq \varepsilon e^{K\tau^{\frac{1}{2}}M^{\frac{1}{2}}}.  & \nonumber
\end{align}

Let $\varepsilon\rightarrow 0$, then $\Theta_{n_k}(s;0,x)\xrightarrow{k}\Theta(s;0,x)$. Then,
\begin{align}\label{tempory}
\|\Theta_{n_k}\Theta_{n_k}^T- \Theta\Theta^T\| &\leq \|\Theta_{n_k}\Theta_{n_k}^T- \Theta\Theta_{n_k}^T\|+\|\Theta\Theta_{n_k}^T-\Theta\Theta^T\| &\nonumber \\
&\leq \|\Theta_{n_k}- \Theta\|\|\Theta_{n_k}^T\|+\|\Theta\|\|\Theta_{n_k}^T-\Theta^T\| &\nonumber\\
&\leq \tilde{M}\|\Theta_{n_k}- \Theta\|+\tilde{M}\|\Theta_{n_k}^T-\Theta^T\|\xrightarrow{k} 0, &
\end{align}
since $\Theta_{n_k}\xrightarrow{k}\Theta$ and $\Theta_{n_k}^T\xrightarrow{k}\Theta^T$.

{\bf Step 4.} We claim that $R_{n_k}\xrightarrow{k}R$.

\vskip1mm Define $A_{n_k}=\Theta_{n_k}(\tau;0,x)\Theta_{n_k}^T(\tau;0,x)=\left(
                                                        \begin{array}{c}
                                                          a^{n_k}_{ij}(\tau;0,x) \\
                                                        \end{array}
                                                      \right)_{3\times 3},
                                                      A=\Theta(\tau;0,x)\Theta^T(\tau;0,x)=\left(
                                                        \begin{array}{c}
                                                          a_{ij}(\tau;0,x) \\
                                                        \end{array}
                                                      \right)_{3\times 3}
$.
\vskip1mm By (\ref{tempory}), we know $a^{n_k}_{ij}(\tau;0,x)\xrightarrow{k} a_{ij}(\tau;0,x)$ for $i,j=1,2,3$. Now, we simply denote  $a^{n_k}_{ij}(\tau;0,x)$ and  $a_{ij}(\tau;0,x)$ by $a^{n_k}_{ij}$ and $a_{ij}$, respectively.
\vskip1mm By Lemma \ref{lemm3-5}, we obtain that
\begin{align}\label{eq-Fi2}
&{\rm det}(\Theta_{n_k}(\tau;0,x))=e^{-\int_0^\tau \sum\limits_{i=1}^3(v_{n_k}(\eta_{n_k}(s;0,x),s))_{i,(\eta_{n_k})_i}ds}\neq 0. &
\end{align}
\vskip1mm By (\ref{eq-Fi2}), we know that $\Theta_{n_k}\Theta_{n_k}^T,\Theta\Theta^T\in SPD(3)$. Let $\lambda_{n_k}^{(1)}\geq \lambda_{n_k}^{(2)}\geq \lambda_{n_k}^{(3)}>0$, $\lambda^{(1)}\geq \lambda^{(2)}\geq \lambda^{(3)}>0$ be the eigenvalues of $\Theta_{n_k}\Theta_{n_k}^T,\Theta\Theta^T$, respectively.
\vskip1mm What's more, by (\ref{eq-Fi2}), we obtain that
\begin{align}\label{cvb}
\det(A_{n_k})=&e^{-2\int_0^\tau \sum\limits_{i=1}^3(v_{n_k}(\eta_{n_k}(s;t,x),s))_{i,(\eta_{n_k})_i}ds}\leq e^{2\int_0^\tau\|\nabla v(\cdot,s)\|ds}
\leq e^{2\int_0^\tau\|v_{n_k}(\cdot,s)\|_{[C^{1}(\Omega)]^3}ds}\nonumber \\
\leq & e^{2C\int_0^\tau\|v_{n_k}(\cdot,s)\|_{[F^\alpha_{L,0}(\Omega)]^3}ds}
\leq e^{2K\int_0^\tau\|\nabla^\alpha v_{n_k}(\cdot,s)\|_{L^2(\Omega)}ds}
\leq e^{2K\tau^{\frac{1}{2}}[\int_0^\tau\|\nabla^\alpha v_{n_k}(\cdot,s)\|^2_{L^2(\Omega)}ds]^{\frac{1}{2}}}\nonumber \\
\leq & e^{2K\tau^{\frac{1}{2}}M^{\frac{1}{2}}}\triangleq M_1<+\infty.&
\end{align}
\vskip1mm In a similar way, we can obtain that
\begin{align}\label{cvb1}
\det(A_{n_k})\geq \frac{1}{M_1}.
\end{align}
\vskip1mm By singularity decomposition theorem\cite{JOR}, we can find two $3\times 3$ orthogonal matrix $U_{n_k}$, $V_{n_k}$ such that
$\Theta_{n_k}=U_{n_k}S_{n_k}V_{n_k}^T$, where $S_{n_k}=\mathrm{diag}\left(\sqrt{\lambda_{n_k}^{(1)}},\sqrt{\lambda_{n_k}^{(2)}},\sqrt{\lambda_{n_k}^{(3)}}\right)$, the columns of $U_{n_k}$, $V_{n_k}$ are orthogonal eigenvectors of $\Theta_{n_k}\Theta_{n_k}^T$ and $\Theta_{n_k}^T\Theta_{n_k}$, respectively.

Similarly, $\Theta=USV^T$, where $S=\mathrm{diag}\left(\sqrt{\lambda^{(1)}},\sqrt{\lambda^{(2)}},\sqrt{\lambda^{(3)}}\right)$, the columns of $U$, $V$ are orthogonal eigenvectors of $\Theta\Theta^T$ and $\Theta^T\Theta$, respectively.

\vskip1mm Then, $A_{n_k}=\Theta_{n_k}\Theta_{n_k}^T=U_{n_k}S_{n_k}^2U_{n_k}^T$ , $A_{n_k}^{-1}=U_{n_k}S_{n_k}^{-2}U_{n_k}^T$ and $A=\Theta\Theta^T=US^2U^T$, $A^{-1}=US^{-2}U^T$. Hence,

\begin{align}
\|A_{n_k}^{-1}\|\leq&\|U_{n_k}\|\|S_{n_k}^{-2}\|\|U_{n_k}^T\|
\leq \|U_{n_k}\|^2\left[\frac{1}{\lambda^{(1)}_{n_k}}+\frac{1}{\lambda^{(2)}_{n_k}}+\frac{1}{\lambda^{(3)}_{n_k}}\right]&\nonumber \\
\leq&\|U_{n_k}\|^2\frac{\lambda^{(1)}_{n_k}\lambda^{(2)}_{n_k}+\lambda^{(1)}_{n_k}\lambda^{(3)}_{n_k}+\lambda^{(2)}_{n_k}\lambda^{(3)}_{n_k}}{\lambda^{(1)}_{n_k}\lambda^{(2)}_{n_k}\lambda^{(3)}_{n_k}}
=\|U_{n_k}\|^2\frac{\lambda^{(1)}_{n_k}\lambda^{(2)}_{n_k}+\lambda^{(1)}_{n_k}\lambda^{(3)}_{n_k}+\lambda^{(2)}_{n_k}\lambda^{(3)}_{n_k}}{\det(A_{n_k})}&\nonumber
\end{align}
\begin{align}\label{ghjk}
\leq&\|U_{n_k}\|^2\frac{[\lambda^{(1)}_{n_k}+\lambda^{(2)}_{n_k}+\lambda^{(3)}_{n_k}]^2}{\det(A_{n_k})}
=\|U_{n_k}\|^2\frac{[\mathrm{tr}(A_{n_k})]^2}{\det(A_{n_k})}
\leq \|U_{n_k}\|^2\frac{\|A_{n_k}\|^2}{\det(A_{n_k})}&\nonumber \\
\leq&27\tilde{M}^2M_1\triangleq M_2<+\infty,&
\end{align}
by (\ref{cvb}), (\ref{cvb1}) and {\bf Step 1}, where $\mathrm{tr}(A)$ denote the trace of matrix $A$. Note that here we use the equalities
\begin{align}
\lambda^{(1)}_{n_k}+\lambda^{(2)}_{n_k}+\lambda^{(3)}_{n_k}\equiv\mathrm{tr}(A_{n_k}), \ \lambda^{(1)}_{n_k}\lambda^{(2)}_{n_k}\lambda^{(3)}_{n_k}\equiv\det(A_{n_k}).
\end{align}
\vskip1mm Similarly , we know that $\|A^{-1}\|\leq M_2$.
\vskip1mm By (\ref{ghjk}), we obtain that
\begin{align}\label{nerf}
\|A_{n_k}^{-1}-A^{-1}\|=&\|A_{n_k}^{-1}(A-A_{n_k})A^{-1}\|\leq \|A_{n_k}^{-1}\|\|A-A_{n_k}\|\|A^{-1}\|\xrightarrow{k} 0,&
\end{align}
since $A-A_{n_k}=\Theta_{n_k}\Theta_{n_k}^T-\Theta\Theta^T\xrightarrow{k} 0$ by (\ref{tempory}).

\vskip1mm Hence, $A_{n_k}^{-1}\xrightarrow{k} A^{-1}$.
\vskip1mm Since $A_{n_k}^{-1}=U_{n_k}S_{n_k}^{-2}U_{n_k}^T=\left[U_{n_k}S_{n_k}^{-1}U_{n_k}^T\right]\left[U_{n_k}S_{n_k}^{-1}U_{n_k}^T\right]\triangleq B_{n_k}B_{n_k}$, $A^{-1}=US^{-2}U^T=\left[US^{-1}U^T\right]\left[US^{-1}U^T\right]\triangleq BB$, then $B_{n_k},B\in SPD(3)$. Then, the Minkowskii inequality[15, Equation (1.1)] implies
\begin{align}
\left[\det(B_{n_k}+B)\right]^{\frac{1}{3}}\geq \left[\det(B_{n_k})\right]^{\frac{1}{3}}+\left[\det(B)\right]^{\frac{1}{3}}\geq \left[\det(B)\right]^{\frac{1}{3}}.
\end{align}

This is, $\det(B_{n_k}+B)\geq \det(B)=\frac{1}{\lambda^{(1)}\lambda^{(2)}\lambda^{(3)}}>0$.
\vskip1mm Further more, we have
\begin{align}\label{nerfhj1}
\|B_{n_k}+B\|\leq &\|B_{n_k}\|+\|B\|
\leq\|U_{n_k}\|^2\left[\frac{1}{\sqrt{\lambda_{n_k}^{(1)}}}+\frac{1}{\sqrt{\lambda_{n_k}^{(2)}}}+\frac{1}{\sqrt{\lambda_{n_k}^{(3)}}}\right]\nonumber&\\
&+\|U\|^2\left[\frac{1}{\sqrt{\lambda^{(1)}}}+\frac{1}{\sqrt{\lambda^{(2)}}}+\frac{1}{\sqrt{\lambda^{(3)}}}\right]&\nonumber\\
\leq&\|U_{n_k}\|^2\frac{\lambda_{n_k}^{(1)}+\lambda_{n_k}^{(2)}+\lambda_{n_k}^{(3)}}{\sqrt{\lambda_{n_k}^{(1)}\lambda_{n_k}^{(2)}\lambda_{n_k}^{(3)}}}+\|U\|^2\frac{\lambda^{(1)}+\lambda^{(2)}+\lambda^{(3)}}{\sqrt{\lambda^{(1)}\lambda^{(2)}\lambda^{(3)}}}\nonumber&\\
=&\|U_{n_k}\|^2\frac{\mathrm{tr}(A_{n_k})}{\sqrt{\det(A_{n_k})}}+\|U\|^2\frac{\mathrm{tr}(A)}{\sqrt{\det(A)}}&\nonumber \\
\leq& 9M_2\sqrt{M_1}+9M_2\sqrt{M_1}\leq 18M_2\sqrt{M_1}\triangleq M_3<+\infty.&
\end{align}
\vskip1mm By (\ref{nerf}), we obtain that
\begin{align}\label{nerf1}
A_{n_k}^{-1}-A^{-1}=B_{n_k}^2-B^2=(B_{n_k}+B)(B_{n_k}-B)\xrightarrow{k} 0.
\end{align}
\vskip1mm By (\ref{nerf1}) and Lemma 2.5 in \cite{HH}, we obtain that $(\Theta_{n_k}\Theta_{n_k}^T)^{-\frac{1}{2}}=B_{n_k}\xrightarrow{k} B=(\Theta\Theta^{T})^{-\frac{1}{2}}$.
\vskip1mm This yields
\begin{align}\label{nerf2}
&\|\Theta_{n_k}(\Theta_{n_k}\Theta_{n_k}^T)^{-\frac{1}{2}}-\Theta(\Theta\Theta^T)^{-\frac{1}{2}}\|&\nonumber \\
\leq&\|\Theta_{n_k}\|\|(\Theta_{n_k}\Theta_{n_k}^T)^{-\frac{1}{2}}-(\Theta\Theta^T)^{-\frac{1}{2}}\|+\|\Theta_{n_k}-\Theta\|\|(\Theta\Theta^T)^{-\frac{1}{2}}\|\xrightarrow{k} 0.&
\end{align}
\vskip1mm So, $R_{n_k}\xrightarrow{k} R$.

{\bf Step 5.} We claim that $\|T\diamond h_{n_k}(\cdot)-D(\cdot)\|^2_{L^2(\Omega)}\xrightarrow{k} \|T\diamond h(\cdot)-D(\cdot)\|^2_{L^2(\Omega)}$.

Let $x\in \Omega\setminus h^{-1}(\Delta_T)$, by (\ref{eq3-24}) and {\bf Step 4}, it yields
\begin{align}\label{final}
\|T\diamond h_{n_k}(x)-T\diamond h(x)\|&=\|R_{n_k} [T\circ h_{n_k}(x)] R_{n_k}^T-R [T\circ h(x)] R^{ T}\| \nonumber &\\
&\leq\|R_{n_k} [T\circ h_{n_k}(x)] R_{n_k}^T-R [T\circ h_{n_k}(x)] R^{T}_{n_k }\| \nonumber &\\
&+\|R [T\circ h_{n_k}(x)] R^{T}_{n_k }-R [T\circ h(x)] R^{T}_{n_k }\| \nonumber &\\
&+\|R [T\circ h(x)] R^{T}_{n_k }-R [T\circ h(x)] R^{ T}\| \nonumber &\\
&\leq\|R_{n_k}-R\|\| T\circ h_{n_k}(x)\|\| R^{T}_{n_k }\| \nonumber &\\
&+\|R\|\| T\circ h_{n_k}(x) - T\circ h(x)\|\|R^{T}_{n_k }\| \nonumber &\\
&+\|R\|\|T\circ h(x)\|\| R^{ T}_{n_k}- R^{ T}\|\xrightarrow{k} 0. &
\end{align}

This implies, $\|T\diamond h_{n_k}(x)-D(x)\|^2\xrightarrow{k} \|T\diamond h(x)-D(x)\|^2$ for all $x\in\Omega\setminus h^{-1}(\Delta_T)$.

On the other hand, $\|T\diamond h_{n_k}(\cdot)-D(\cdot)\|^2\leq\mathop {\max }\limits_{x\in \Omega} \|T(x)-D(x)\|^2=G\in L^1(\Omega)$. Besides, by {\bf (ii)} in Lemma \ref{lemm3-5}, $h^{-1}(\Delta_T)$ is a set of measure 0. Now by Dominance Theorem[5, Theorem 5 in Appendix E], it yields
\begin{align}\label{dom}
\|T\diamond h_{n_k}(\cdot)-D(\cdot)\|_{L^2(\Omega)}^2\xrightarrow{k} \|T\diamond h(\cdot)-D(\cdot)\|_{L^2(\Omega)}^2.
\end{align}

It follows from (\ref{{neweq3-11}}) and (\ref{dom}) that, $H(v)$ is a l.w.c functional. That is,
\begin{align}
\lim_{n_k\rightarrow \infty}\inf H(v_{n_k})\geq H(v).
\end{align}

This concludes the claim (\ref{claim}).

By Lemma \ref{mainlemm}, there exists a global minimizer $\bar{{v}}(x,s)\in B_M$ such that $H(\bar{{v}})=\inf\limits_{v\in B_M} H(v)=\inf\limits_{v\in\mathcal{F}} H(v)=\min\limits_{v\in\mathcal{F}} H(v)$. That is, $\bar{{v}}(x,s)$ is a solution of (\ref{eq1-4}).

\vskip1mm For the above minimizer $\bar{{v}}(x,s)\in\mathcal{F}$, by Lemma \ref{lemm3-4}, we know that there exists a unique $\bar{{\eta}}(s;t,x)\in C([0,\tau],\bar{\Omega})$ such that
\begin{align}
\frac{d\bar{{\eta}}(s;t,x)}{ds}=\bar{{v}}(\bar{{\eta}}(s;t,x),s), \ \ \ \ \bar{{\eta}}(t;t,x)=x.
\end{align}
\vskip1mm Furthermore, by Lemma \ref{lemm3-4} and (\ref{eq1-3}), we know the mapping $\bar{{h}}(x)=\bar{\eta}(0;\tau,x)\in [C^{[\alpha-1.5],\lambda}(\Omega)]^3$ with $\nabla_x\bar{{h}}(x)$ given by Lemma \ref{lemm3-5}, where $0<\lambda<\alpha-[\alpha]$. Moreover, by {\bf(ii)} in Lemma \ref{lemm3-5}, we know $\bar{{h}}$ is a $1$-to-$1$ and onto mapping.\qed

\vskip1mm  {\bf Acknowledgements.}{\ \ This paper is partially based on some results of my P.H.D dissertation supervised by Professor Huan-Song Zhou in Wuhan Institute of Physics and Mathematics, Chinese Academy of Sciences.}\\

%\begin{remark}
%By Lemma \ref{th-jk6}, we know $det(J_x)>0$. That implies $\lambda_i>0\ (i=1,2,\cdots,n)$.
%\end{remark}
%\begin{theorem}\label{th-APPen}
%Consider the singularity decomposition,$\mathbf{M}=\mathbf{U}\mathbf{S}\mathbf{V}^T$, for any square or tall-rectangular matrix, %i.e,$\mathbf{M}\in \mathbb{R} ^{n\times k}$ with $n\geq k$
%\vskip1mm (1) The singular values are unique and for distinct positive singular values $s_j>0$, the $j^{th}$ columns of %$\mathbf{U},\mathbf{V}$ are also unique up to a sign change of both columns.
%\vskip1mm (2) For any repeated and positive singular values, say $s_i=s_{i+1}=\cdots=s_j>0$ are all singular values equals to $s_i$, %the corresponding columns of $\mathbf{U},\mathbf{V}$ are unique up to any rotation/reflection applied to both sets of columns(i.e %$\mathbf{U}_{*,i:j}\rightarrow \mathbf{U}_{*,i:j}\mathbf{W},\mathbf{V}_{*,i:j}\rightarrow \mathbf{V}_{*,i:j}\mathbf{W}$ for some %orthogonal matrix $\mathbf{W}$).
%\vskip1mm \ \ \ \ \ \ \ \ \ \ \ \ \ \ \ \ \ \ \ \ \ \ \ \ \ \ \ \ \ \ \ \ \ \ \ \ \ \ \ \ \ \ \ \ \ \ \ \ \ \ \ \ \ \ \ \ \ \ \ \ \ \ \ \ \ \ \ \ \ \ \ \ \ \ \ \ \ \ \ \ \ \ \ \ \ \ \ \ \ \ \ \ \ \ \ \ \ \ \ \ \ \ \ \ \ \ \ \ \ \ \ \ \ \ \ \ \ \ \ \ \ \mbox{$\blacksquare$}
%\end{theorem}

\end{document}